\documentclass[11pt,tgrind,psbox]{article}
\usepackage{latexsym}
\usepackage{graphics}
\usepackage{amsmath}
\usepackage{xspace}
\usepackage{amssymb}
\usepackage{psfrag}
\usepackage{epsfig}
\usepackage{fullpage}
\usepackage{amsmath,amsthm}
\usepackage{epsfig}
\usepackage{pst-all}
\usepackage{bm}

\usepackage{times} \usepackage{amssymb} \usepackage{amsfonts} \usepackage{amsmath} \usepackage{latexsym}

\newcounter{enumerate_t} 
\newcommand{\ignore}[1]{}     \makeatletter
\renewcommand{\section}{\@startsection{section}{1}{0pt}{-8pt}{6pt}{\large\bf}}
\renewcommand{\subsection}{\@startsection{subsection}{2}{0pt}{-6pt}{-5pt}{\normalsize\bf}}
\renewcommand{\subsubsection}{\@startsection{subsubsection}{3}{0pt}{-12pt}{-5pt}{\normalsize\bf}} \makeatother  
\newenvironment{Proof}{\noindent{\sc Proof. }}{$\blacksquare$\medskip}
\newtheorem{claim}{Claim} 
\newtheorem{lemma}[claim]{Lemma}
 \newtheorem{theorem}[claim]{Theorem} 
 \newtheorem{corollary}[claim]{Corollary}

\newcommand{\reals}{{\mathbb R}}    

\newcommand{\p}{{\mathbb P}}

\newcommand{\E}{\mathbb{E}}

\def\qquad{\hskip2em\relax}
\DeclareMathOperator*{\argmin}{arg\,min}

\usepackage[colorlinks=true,breaklinks=true,bookmarks=true,urlcolor=blue,
     citecolor=blue,linkcolor=blue,bookmarksopen=false,draft=false]{hyperref}

\begin{document}




\title{Asymptotic Optimality of Constant-Order Policies for Lost Sales Inventory Models with Large Lead Times}


\author{
{\sf David A. Goldberg}
\thanks{H. Milton Stewart School of Industrial and Systems Engineering, Georgia Institute of Technology, e-mail: {\tt dgoldberg9@isye.gatech.edu}}\ \ \ \ \ 
{\sf Dmitriy A. Katz-Rogozhnikov}
\thanks{Business Analytics and Mathematical Sciences, IBM T.J. Watson Research Center, e-mail: {\tt dkatzrog@us.ibm.com}}\ \ \ \ \ 
{\sf Yingdong Lu}
\thanks{Business Analytics and Mathematical Sciences, IBM T.J. Watson Research Center, e-mail: {\tt yingdong@us.ibm.com}}\\
{\sf Mayank Sharma}
\thanks{Business Analytics and Mathematical Sciences, IBM T.J. Watson Research Center, e-mail: {\tt mxsharma@us.ibm.com}}\ \ \ \ \ 
{\sf Mark S. Squillante}
\thanks{Business Analytics and Mathematical Sciences, IBM T.J. Watson Research Center, e-mail: {\tt mss@us.ibm.com}}
}

\maketitle


\begin{abstract}
Lost sales inventory models with large lead times, which arise in many
practical settings, are notoriously difficult to optimize due to the curse of
dimensionality.
In this paper we show that when lead times are large, a very simple
constant-order policy, first studied by Reiman \cite{Reiman04}, performs nearly optimally.
The main insight of our work is that when the lead time is very large, such a
significant amount of randomness is injected into the system between when an
order for more inventory is placed and when the order is received, that
``being smart'' algorithmically provides almost no benefit.
Our main proof technique combines a novel coupling for suprema of random walks
with arguments from queueing theory.
\end{abstract}



\maketitle

%


%
%
%

\section{Introduction}
\label{sec:intro}
In this paper we consider a stochastic inventory control problem under the
so-called single-item, periodic-review, lost-sales model with positive lead
times and independent and identically distributed (i.i.d.) demand.
This model is based on sales being lost whenever there is insufficient supply
to fulfill demand, i.e., unfulfilled demand is lost rather than being carried
over, or backlogged, to a later time.
Furthermore, there is a constant delay of $L > 0$ periods (i.e., a single
lead time) between when an order for additional inventory is placed and when
that inventory is received.
The problem then is to determine the best policy for a series of orders
across a planning horizon comprised of a finite number of discrete time
periods, with the goal of minimizing cost in expectation.

The cost structure of this model consists of a per-unit penalty for lost sales
due to unfulfilled demand within each period and a per-unit cost for holding
excess inventory within each period.
Unlike the corresponding backorder inventory control problem when unfulfilled
demand is fully backlogged from period to period, where the optimal policy
is well known to be an order-up-to policy, the optimal order policy for
the lost-sales inventory model is not known in general, and in fact remains
poorly understood~\cite{Bijvank11}.

Such periodic-review, lost-sales models have a long and rich history in the
operations research, operations management and management science literature.
Here we briefly review some of the most relevant literature, and refer the interested
reader to the recent survey paper of Bijvank \cite{Bijvank11} for a more comprehensive
exposition.
This class of inventory models was first introduced by Bellman \cite{Bellman55}.
Certain properties of the optimal policy were explored for the case of $L = 1$
by Karlin and Scarf \cite{KarlinScarf} and by Yaspan \cite{Yaspan61}, where it
was shown that the  order-up-to policy is not optimal for the lost-sales
inventory model.
Morton \cite{Morton69} extended this analysis to the case of general $L$.
Other properties of the optimal policy, including various notions of convexity
and monotonicity, were explored in \cite{Zipkin08,Zipkin08b,HuhOR10}.
With respect to computation of the optimal policy, the primary approach taken
in the literature is dynamic programming, combined with various heuristics to
speed up computations \cite{Morton71,Zipkin08b}.
However, since the state-space of any such dynamic program grows exponentially
in the lead time, such computations become extremely challenging even for lead times less than
ten~\cite{Zipkin08b}.
Namely, this family of techniques suffers from the curse of dimensionality as
the lead time grows.
Indeed, even for a lead time of four and geometrically distributed demand, Zipkin
\cite{Zipkin08b} reports that computing the optimal policy requires
solving a dynamic program with $228,581$ states.
This is not surprising because the problem at hand and several closely related
problems are known to be NP-complete~\cite{Halman12}.

The difficulty of computing optimal policies for the lost-sales model has led
to a considerable body of work on heuristics.
The computational performance and properties of various algorithms, including
order-up-to policies, have been analyzed by numerous
authors~\cite{Gaver59,Morse59,Morton71,Yaspan72,Pressman77,Nahmias79,Downs01,Johansen01,Janakiraman04,Johansen08,Bijvank10}.
With respect to policies that have provable performance guarantees, the
breakthrough work of Levi et al. \cite{Levi08} proved that a certain
dual-balancing heuristic, inspired by previous results for other
models~\cite{Levi08c, Levi08b}, yields a policy whose cost is always within
a factor of two of optimal.
Huh et al. \cite{huh2009asymptotic} show that in a certain scaling regime, in which
the ratio of the lost-sales penalty to the holding cost asymptotically tends
to infinity, an order-up-to policy is asymptotically optimal;
and a similar result has been recently derived by Lu et al. \cite{LuSqYa12}.
Using a very different approach, Halman et al.\ provide an approximate dynamic
programming algorithm that, combined with ideas from discrete convexity,
yields a so-called fully polynomial-time approximation scheme for
various related inventory control problems \cite{Halman09,Halman12}.
These techniques were recently extended to lost-sales models with positive
lead times (as considered in this paper) by Chen et al. \cite{Chen12a}, who
provide a pseudo-polynomial-time additive approximation algorithm.
Namely, under a suitable encoding scheme, an algorithm is presented that,
for any $\epsilon > 0$, returns a policy whose performance differs additively
from that of the optimal policy by at most $\epsilon$, in time which is
polynomial in $\epsilon^{-1}$ if the overall encoding length of the problem
is held fixed while $\epsilon$ is varied, and otherwise is pseudo-polynomial
in the overall encoding length (which grows with the lead time $L$);
we refer the reader to \cite{Chen12a} for details.
In a follow-up study~\cite{Chen12b}, the authors prove several interesting
integrality results for these and related models.

The work closest to our own is that of Reiman \cite{Reiman04}, who studies a very
simple policy for a certain continuous-review, lost-sales model with positive
lead times and demand arriving as a Poisson process.
In particular, the author analyzes an open-loop constant-order policy, which
at time 0 selects an interval size $\tau$ and simply orders a single unit of
inventory every $\tau$ time units.
The author observes that this simple policy can be analyzed as a $D/M/1$ queue,
and goes on to perform an interesting asymptotic analysis, showing that for any
fixed holding cost and lost-sales penalty, there exists a critical lead-time
value $L^*$ such that (s.t.): (i) for all lead times less than $L^*$, the best
base-stock policy outperforms the best constant-order policy; and (ii) for all
lead times greater than $L^*$, the best constant-order policy outperforms the
best base-stock policy.
The author makes no attempt to compare either policy to the true optimal
policy, which he notes is unknown.

Of course, there is no a priori reason to believe that such a simple
constant-order policy should be nearly optimal.
However, numerical results from a recent study by Zipkin \cite{Zipkin08b}, in which
the optimal policy is computed for a lost-sales model with i.i.d.\ demand and positive lead times (nearly identical to
the model we consider, but with discounting), show that the constant-order
policy (in which the same fixed constant is ordered in every time period) can
perform surprisingly well.
More precisely, in numerical experiments for a lead time of four, the
constant-order policy always incurs an expected cost at most twice that
incurred by the optimal policy; in 62.5\% of the cases, the constant-order
policy incurs a cost at most $1.33$ times that incurred by the optimal policy;
and in 38\% of the cases, it incurs a cost at most $1.12$ times that incurred
by the optimal policy.
This begs the question of how such a simple policy can perform so well on
reasonable problem instances.

In the present paper we derive theoretical results that shed light on this
and related phenomena.
Specifically, we prove that, as the lead time grows (with the demand
distribution, lost-sales penalty and holding cost remaining fixed),
the best constant-order policy is in fact asymptotically optimal.
We also establish explicit bounds on how large the lead time should be to
ensure that the best constant-order policy incurs an expected cost of at
most $1+\epsilon$ times that incurred by the optimal policy.
To the best of our knowledge, this is the first algorithm proven to be within
$1 + \epsilon$ of optimal for lost-sales models when the lead time is large,
and whose runtime does not grow with the lead time.
The main insight of our work is that when the lead time is very large, such a
significant amount of randomness is injected into the system between when an
order for more inventory is placed and when the order is received, that
``being smart'' algorithmically provides almost no benefit.
Our main proof technique combines a novel coupling for suprema of random walks
with arguments from queueing theory.
Since this simple policy succeeds exactly when known algorithms start running
into trouble due to the curse of dimensionality, our results open the door for
the creation of ``hybrid'' algorithms that use more elaborate forms of dynamic
programming when the lead time is small, and gradually transition to less
computationally intensive algorithms (with the constant-order policy at the
extreme) as the lead time grows.

\subsection*{Outline of paper and overview of proof.}
The remainder of this paper, including the underlying proofs of our main results, is organized as follows.
Section~\ref{model} formally defines the model of study and Section~\ref{main} states our main results (Theorem~\ref{mainresult} and Corollary~\ref{maincorr}),
namely the asymptotic optimality of the constant-order policy and associated explicit performance guarantees.
In Section~\ref{Constant-order policy}, we explicitly describe the dynamics and associated costs for a general policy over any consecutive $L$ time periods,
if one conditions on the pipeline-vector and inventory at the start of those $L$ time periods, in terms of the maxima of various partial sums (Lemma~\ref{explic1}).
We then customize this result to the constant-order policy (Lemma~\ref{constcost}), under which various associated expressions simplify considerably.

In Section~\ref{lowerbound}, we develop lower bounds on the cost incurred by a particular optimal policy $\overline{\pi}$ described in \cite{Zipkin08},
which never orders more than a certain quantity that depends on the underlying costs and demand distribution, but not on the leadtime $L$ nor planning horizon $T$.
First, we formulate a lower bound on the cost incurred by $\overline{\pi}$ over any consecutive $L$ time periods by supposing that $\overline{\pi}$ was able to choose the state of the system at the start of those $L$ time periods to be as favorable as possible.
This results in a ``best-case'' pipeline vector $\mathbf{x}^*$ and inventory level ${\mathcal I}^*$, which can be described as the solution to an appropriate optimization problem (\ref{getlow1}), and then can be used together with Lemma~\ref{explic1} to derive a lower bound for the cost incurred by $\overline{\pi}$ over any consecutive $L$ periods (Lemma~\ref{lowerblemma}), again in terms of the maxima of various partial sums.

The second lower bound formulated in Section~\ref{lowerbound} (Lemma~\ref{lowerblemma2}) represents the most critical step of the entire proof,
where we reason as follows.
Ultimately, we wish to show that the performance of an appropriate constant-order policy nearly matches the lower bound established in Lemma~\ref{lowerblemma}.
To accomplish this, we first note that if we wanted to select a constant-order policy which came close to matching the aforementioned lower bound, a natural approach would be to select the constant-order policy that ``best mimics'' the pipeline vector $\mathbf{x}^*$, i.e., by defining $r^* \stackrel{\Delta}{=} L^{-1} \sum_{i=1}^L x^*_i$ and considering the policy that orders $r^*$ in every period (assuming $r^* < \E[D]$).
Second, we note that if one compares the maxima appearing in the lower bound of Lemma~\ref{lowerblemma}, and the maxima appearing in the dynamics of the constant-order policy when $r^*$ is ordered in every period (as described in Lemma~\ref{constcost}), the associated expressions are markedly similar.
A fundamental difficulty in precisely comparing these expressions is that the indices at which various maxima appearing in Lemma~\ref{lowerblemma} attain their suprema may depend on fluctuations in the vector $\mathbf{x}^*$, while the maxima appearing in Lemma~\ref{constcost} do not have this property, as the associated pipeline vector is constant.
To remedy this, we describe an explicit coupling between the maxima that appear in Lemma~\ref{lowerblemma} and the maxima that appear in Lemma~\ref{constcost} when $r^*$ is ordered in every period.
Indeed, instead of computing the ``true values'' of the maxima appearing in Lemma~\ref{lowerblemma}, we derive bounds by examining the associated expressions not at the index at which they attain their supremum, but at a different index, namely the index at which a corresponding expression appearing in Lemma~\ref{constcost} attains its supremum.
In this way, we are able to derive a second lower bound (Lemma~\ref{lowerblemma2}), which is much closer to the corresponding expression appearing in Lemma~\ref{constcost}, and highly amenable to analysis.

Then in Section~\ref{boundthedifference} we use Lemma~\ref{lowerblemma2}, along with a careful analysis of various quantities associated with certain random walks and their suprema (Lemma~\ref{boundindex}), to explicitly bound the performance gap between the aforementioned constant-order policy and $\overline{\pi}$, over any consecutive $L$ periods (Theorem~\ref{afewgoodbounds}).
In Section~\ref{provesec}, we use Theorem~\ref{afewgoodbounds} to complete the proof of our main results.
Along these lines, we first prove that $r^*$ is bounded away from $\E[D]$ (Lemma~\ref{boundr}), as the bounds of Theorem~\ref{afewgoodbounds} depend sensitively on this gap.
The associated proof proceeds by contradiction in which it is shown that if $r^*$ is ``too close'' to $\E[D]$, then a certain ``sharpening'' of the central limit theorem, known as Stein's lemma, guarantees that too large of a holding cost is incurred.
Finally, we note that every constant-order policy is dominated by a certain ``best-possible'' constant-order policy.
Combining all of the above with a straightforward asymptotic analysis completes the proof.

Closing remarks and directions for future research are presented in Section~\ref{conclusion}.
We also include a short technical appendix in Section~\ref{appsec}.

\section{Model description and problem statement}\label{model}
Let us consider a standard lost-sales inventory optimization problem.
One is given as input the holding cost, lost-demand penalty, planning horizon,
lead time, and demand distribution.
The problem then is to determine the optimal ordering policy to control
inventory in the so-called single-item, periodic-review, lost-sales model.

Specifically, we consider the following lost-sales model and associated
optimization problem.
Let $h$ denote the per-unit holding cost and $c$ the per-unit lost-demand
penalty, where we assume $c,h > 0$.
Time is slotted, with the planning horizon and lead time respectively comprised of $T$ and $L$ periods s.t.\ $T > L > 0$.
The demand in each period $t$, denoted by $D_t$, is assumed to be i.i.d.\ and governed by a non-negative demand distribution ${\mathcal D}$ with finite third moment.  To prevent certain degenerate situations, we further assume that ${\mathcal D}$ is not deterministic.

At the start of each time period $t$ there is an amount of inventory $I_t$ available.
There is also an $L$-dimensional vector of pipeline inventory $\mathbf{x}_t = (x_{1,t}, x_{2,t}, \ldots, x_{L,t})$
that captures the orders placed before period $t$ but not yet received prior to period $t$.
The system dynamics for time period $t$ then proceed as follows.
First, a new amount $x_{1,t}$ of goods is added to the inventory.
Second, before seeing the demand of period $t$, an order for more inventory is placed.
After placing this order, the pipeline inventory vector is updated in a manner analogous to that of a queue:
the front entry $x_{1,t}$ is removed, all other entries move up one position (i.e., $x_{i, t+1} = x_{i+1, t}$ for $i=1, 2, \ldots, L-1$),
and the new order is appended at the end, becoming $x_{L,t+1}$.

The order placed at time $t$ must be a function (albeit possibly a random function) only of the planning horizon $T$, the current time $t$,
the inventory level $I_t$ at the start of period $t$,
the pipeline vector $\mathbf{x}_t$ at the start of period $t$,
and the model primitives $L, h, c, {\mathcal D}$.
In particular, the ordering decision at time $t$ \emph{cannot} depend on the
realizations of future demand.
We call all such policies {\it admissible} policies, and denote the family of
admissible policies by $\Pi$.

Next, a random demand $D_t$ is realized from ${\mathcal D}$.
At the end of time period $t$ (but before the start of period $t+1$), the inventory or lost-sales costs are incurred as follows.
The amount of excess inventory at the end of period $t$ is given by $I_{t+1} = ( I_t + x_{1,t} - D_t )^+$, noting that $D_t$ is independent of $I_t + x_{1,t}$.
Conversely, the amount of lost demand (due to not having enough inventory on hand) in period $t$ is denoted by $N_t \stackrel{\Delta}{=} ( I_t + x_{1,t} - D_t )^-$.
Then, the holding cost incurred (for storing excess inventory) equals $h I_{t+1}$ and the lost-sales penalty incurred (for lost demand) equals $c N_t$,
noting that at most one of $I_{t+1}$ and $N_t$ is positive for any period $t$.


The goal of the inventory planner is to minimize the expected cost incurred over the entire planning horizon.  Let $\mathbf{1}$ denote the vector of all ones and $\mathbf{0}$ the vector of all zeros, where the dimension is to be inferred from context.
Let us suppose $I_0 = 0$ and $\mathbf{x}_0 = \mathbf{0}$, which should be assumed throughout as the given initial conditions, and let us further define
$$
C_t \;\; \stackrel{\Delta}{=} \;\; h I_{t+1} + c N_t \;\; = \;\; h  ( I_t + x_{1,t} - D_t )^+ + c  ( I_t + x_{1,t} - D_t )^-.
$$
The planner then wishes to find the policy $\pi \in \Pi$ that minimizes $\E[\sum_{t = L + 1}^{T + L} C_t] $,
where the expectation is over the random demand and any random decisions taken by policy $\pi$, and
where we suppose that ordering decisions are made only in periods $1,\ldots,T$.
Note that we do not penalize a policy for any costs incurred over the first $L$ time periods,
as this cost is completely determined by the initial pipeline vector and random demands.
On the other hand, we do penalize a policy for costs incurred in periods $[T, T+L]$ (recalling $T > L$),
as the cost incurred during these periods is completely determined by the ordering decisions made in periods $1,\ldots,T$ and random demands.
We note that such a convention is consistent with the previous literature; see, e.g., \cite{Zipkin08}.
As a notational convenience and without loss of generality (w.l.o.g.),
we suppose that every policy orders $0$ in periods $T+1,\ldots,T+L$, as these ordering decisions have no effect on the problem's cost.
For a given policy $\pi$, let $\lbrace N^{\pi}_t, C^{\pi}_t, I^{\pi}_t, \mathbf{x}^{\pi}_t ; \, t =1,\ldots,T + L \rbrace$
denote the associated random variables (r.v.s) when policy $\pi$ is implemented
(all constructed on the same probability space).
The corresponding lost-sales inventory optimization problem is then given by
\begin{equation}\label{probstate}
\inf_{\pi \in \Pi} \sum_{t = L + 1}^{T + L} \E[ C^{\pi}_t].
\end{equation}

\section{Main results}
\label{main}
Our main results establish that there exists a very simple constant-order policy which is asymptotically optimal as $L \rightarrow \infty$.
This section formally states these results.

\subsection{Additional definitions and notations.}

Let $D$ denote a r.v.\ governed by ${\mathcal D}$.
Note that if the same deterministic quantity $r < \E[D]$ is ordered in every
period, then the inventory evolves exactly as the waiting time in a $GI/GI/1$
queue (initially empty) with interarrival time distribution ${\mathcal D}$ and
processing time distribution (the constant) $r$;
we refer the reader to \cite{Asmu03} for an excellent discussion of the
dynamics and steady-state properties of the $GI/GI/1$ queue.
Let $I^r_{\infty}$ denote a r.v.\ distributed as the steady-state waiting time in the corresponding $GI/D/1$ queue;
namely, $I^r_{\infty}$ is distributed as $\sup_{k \geq 0}(k r - \sum_{i=1}^k D_i)$.

For two vectors $\mathbf{x},\mathbf{y}$ of equal dimension, we use the notation $\mathbf{x} \leq \mathbf{y}$ to denote component-wise domination, i.e., $x_i \leq y_i$ for all $i$.
Define $Q$ to be the $\frac{c}{c + h}$ quantile of the demand distribution, namely
$$Q \;\; \stackrel{\Delta}{=} \;\; \inf \lbrace s \in {\mathcal R}^+ : \p(D > s) \leq \frac{h}{c + h} \rbrace.$$
We note that $Q$ is the optimal inventory level for the corresponding single-stage newsvendor problem \cite{zipkin2000fundamentals},
i.e., for any policy $\pi$ and any time $t$
$$\E[C^{\pi}_t] \;\; \geq \;\; g \; \triangleq  \; h \E[(Q - D)^+] + c \E[(D - Q)^+].$$

Lastly, let $\sigma$ denote the standard deviation of $D$, and $\zeta \stackrel{\Delta}{=} \E\big[|D - \E[D]|^3\big] \sigma^{-3}$ denote the so-called skewness of $D$.
We then define several functions that will be instrumental for our analysis:
$$m \triangleq \left\lceil \Big(26 \big( 3 \zeta + c (h \sigma)^{-1} \E[D] + 1 \big) \Big)^2 \right\rceil \; , \qquad\qquad z \triangleq \argmin_{v \geq 0} \; \big(h \E\big[I^v_{\infty}\big] - c v\big),$$
and
$$
y (\epsilon) \triangleq \max\left( 2^{14} h (Q + 2^{\frac{3}{2}} \E[D])(\E^2[D] + \E[D^2])^3 \sigma^{-6} m^3 g^{-1} \epsilon^{-1} , \; \Big( 12 c g^{-1} \big( (2 c h^{-1})^{\frac{1}{2}} + 3 \big) \Big)^2 \epsilon^{-2} \right).$$
Although the above quantities are functions of $c,h,{\mathcal D}$, since there
will be no ambiguity, we make this dependence implicit.

\vspace*{0.1in}
\noindent
{\bf Remarks.}
\begin{itemize}
\item
If $\argmin_{v \geq 0} \; \big(h \E\big[I^v_{\infty}\big] - c v\big)$ is not unique, we define $z$ to be the infimum of all such values.
We will later show that $z$ is the best constant possible if the same amount has to be ordered in every time period, in an appropriate sense.
Note that $z \in [0,\E[D])$, since
$\E\big[I^0_{\infty}\big] = 0$ and
$\lim_{r \uparrow \E[D]}\E\big[I^r_{\infty}\big] = \infty$.
\item The function $y(\epsilon)$ captures how large $L$ should be so that our
constant-order policy is within a $(1+\epsilon)$ multiplicative factor of
the optimal policy.
\item 
Note that, for all sufficiently small $\epsilon$, the term
$\left( 12 c g^{-1} \big( (2 c h^{-1})^{\frac{1}{2}} + 3 \big) \right)^2 \epsilon^{-2}$ will dominate
$y(\epsilon)$ due to its quadratic dependence on $\epsilon^{-1}$.
\item We note that our assumption that $D$ has finite third moment is not strictly necessary, but allows for a considerably simplified exposition.  Indeed, in an earlier version of this work \cite{goldberg2012asymptoticv1}, very similar results were presented under the assumption of only finite second moment, but for technical reasons requiring that $D$ had unbounded support.  By combining the arguments of the present paper (which assumes finite third moment but allows $D$ to have bounded support) with the arguments of the aforementioned earlier version, we would in principle obtain asymptotic optimality assuming only a finite second moment.  For simplicity of exposition, we do not formalize this generalization, and leave as an open question the minimal required assumptions on $D$ for such an asymptotic optimality to hold.      
\end{itemize}

\subsection{Formal statement of results.}

For $r \in [0, \E[D])$, let $\pi_{r}$ denote the policy that orders the random amount $I^r_{\infty} + r$
(i.e., the order amount is drawn from the distribution that governs the r.v.\ $I^r_{\infty}$ plus the constant $r$)
in the first time period, and then orders the constant $r$ in all subsequent time periods.
We note that by ordering $I^r_{\infty} + r$ in the first time period, as opposed to $r$, the associated sequence of inventory levels becomes a stationary process, which considerably simplifies our analysis.
With a slight abuse of notation, we still refer to $\pi_r$ as a constant-order policy.
Let $\textrm{OPT}(L,T)$ denote the optimal value of the lost-sales inventory
optimization problem \eqref{probstate} for a given $L$ and $T$.
We then have our main theorem and an important corollary.

\begin{theorem}\label{mainresult}
{\bf} \hspace*{0.1in}
For all $\epsilon \in (0,1)$, $L \geq y(\epsilon)$, and
$T \geq (1 +  3 \epsilon^{-1}) L$,
$$\frac{\E[\sum_{t= L + 1}^{T+ L} C^{\pi_{z}}_t]}{\textrm{OPT}(L,T)} \leq 1 + \epsilon.
$$
\end{theorem}

\begin{corollary}\label{maincorr}
{\bf} \hspace*{0.1in}
$$
\lim_{L \rightarrow \infty} \limsup_{T \rightarrow \infty} \frac{\E[\sum_{t= L + 1}^{T + L} C^{\pi_{z}}_t]}{\textrm{OPT}(L,T)} = 1.
$$
\end{corollary}

In particular, the simple constant-order policy is asymptotically optimal as
$L \rightarrow \infty$.

\section{Inventory dynamics}
\label{Constant-order policy}
In this section, we derive several expressions that explicitly describe the inventory dynamics for any policy over any consecutive $L$ periods,
and then customize these results to the constant-order policy for later use in our proofs.

\subsection{General policy dynamics}\label{gendynsec}
We first explicitly characterize the cost incurred under any given policy during any consecutive $L$ periods.
Although such a characterization is generally well-known (see, e.g., \cite{Zipkin08}), we include a proof for completeness.
Let ${\mathcal I}$ denote the initial inventory level and, for positive integers $j,k$, define $\delta_{j,k}$ to be $1$ if $j = k$ and to be $0$ otherwise.
\begin{lemma}\label{explic1}
{\bf} \hspace*{0.1in}
For any policy $\pi \in \Pi$ and time $\tau \in [1 , T + 1]$,
\begin{eqnarray*}
\E\left[\sum_{t=\tau}^{\tau + L - 1} C^{\pi}_t \Big| \mathbf{x}^{\pi}_{\tau} = \mathbf{x} , I^{\pi}_{\tau} = {\mathcal I} \right] &=&
h \sum_{k=1}^{L} \E\bigg[\max_{j=0,\ldots,k} \bigg( \sum_{i=k+1-j}^k \big( x_i - D_{\tau + i - 1} ) + \delta_{j,k} {\mathcal I} \bigg) \bigg] \\
&& \quad +  \; c \bigg( \E\big[I^{\pi}_{\tau + L} \big| \mathbf{x}^{\pi}_{\tau} = \mathbf{x} , I^{\pi}_{\tau} = {\mathcal I}\big] - {\mathcal I} + L \E[D] - \sum_{i = 1}^{L} x_i \bigg).
\end{eqnarray*}
\end{lemma}

\vspace*{0.1in}
\noindent
{\bf Remark.}
Note that, instead of using the $\delta_{j,k}$ notation to indicate that the initial
inventory is only accounted for in a single term appearing in the associated maximum, we
could have, e.g., considered a transformed set of ``inventory position'' variables, with
each associated variable corresponding to the sum of certain pipeline (and possibly inventory)
variables, as was done in \cite{Zipkin08}.
However, since our subsequent analysis will rely sensitively on the indices at which
various associated maxima are attained, and precisely which pipeline-vector components
appear in the associated partial sums, as well as whether the associated initial inventory
level ${\mathcal I}$ is accounted for in the associated partial sum for which the maximum
is attained (i.e., the maximum occurs at index $k$), we believe that the $\delta_{j,k}$
notation leads to greater clarity of exposition and use this notation throughout.

\vspace*{0.1in}

\begin{Proof}
The result follows from a straightforward, generally well-known (see, e.g., \cite{Zipkin08}), induction that for any $k \in [1 , L]$:
\begin{equation}\label{inventory1}
I^{\pi}_{\tau+k} = \max_{j=0,\ldots,k} \bigg( \sum_{i=k + 1 - j}^k \big( x_i - D_{\tau + i - 1} ) + \delta_{j,k} {\mathcal I} \bigg).
\end{equation}
Note that for any times $t_1 , t_2$ s.t.\ $\tau \leq t_1 < t_2 \leq \tau + L$, the net amount of demand that is met
during $[t_1,t_2 - 1]$ equals $\sum_{t = t_1}^{t_2 - 1} D_t -  \sum_{t = t_1}^{t_2 - 1}  N^{\pi}_t.$
It follows that
\begin{equation}\label{lost2}
\sum_{t = t_1}^{t_2 - 1}  N^{\pi}_t \;\; = \;\; I^{\pi}_{t_2} - I^{\pi}_{t_1} + \sum_{t = t_1}^{t_2 - 1} D_t - \sum_{t = t_1}^{t_2 - 1} x_{t - \tau + 1} .
\end{equation}
Combining (\ref{inventory1}) with (\ref{lost2}) completes the proof.
\end{Proof}

\subsection{Constant-order policy dynamics}\label{constdynsec}
We next customize Lemma~\ref{explic1} to the constant-order policy $\pi_r$.
As previously noted, if the same deterministic quantity $r$ is ordered in every period, the inventory evolves exactly as the waiting time in a $GI/GI/1$ queue with interarrival distribution ${\mathcal D}$ and processing time distribution (the constant) $r$.
Recall that for $r \in [0,\E[D])$, $I^r_{\infty}$ denotes a r.v.\ distributed as the steady-state waiting time in the corresponding stable $GI/D/1$ queue, i.e.,
$I^r_{\infty}$ is distributed as
$\sup_{k \geq 0}(k r - \sum_{i=1}^k D_i)$.
It follows that $\lbrace I^{\pi_r}_{L + k}  ; \, k \geq 1 \rbrace$ is a
stationary sequence of r.v.s, with $I^{\pi_r}_{L + k}$ distributed as
$I^r_{\infty}$ for all $k \geq 1$.

Before explicitly showing how Lemma~\ref{explic1} simplifies when applied to $\pi_r$, it will be useful to introduce some additional notations.
Let $I^r_{1,\infty}$ denote a particular replication of $I^r_{\infty}$
s.t.\ $I^r_{1,\infty}$ and $\lbrace D_i; i \geq 1 \rbrace$ are mutually
independent, where we note that $I^r_{1,\infty}$ will play the role of ${\mathcal I}$ when Lemma~\ref{explic1} is applied to $\pi_r$.
For $k \geq 0$, let us define
$$
W^r_k \stackrel{\Delta}{=} \max_{j=0,\ldots,k} \bigg( j r - \sum_{i=1}^j D_{i} + \delta_{j,k} I^r_{1,\infty} \bigg)
$$
and
$$
i^r_k  \stackrel{\Delta}{=}
\max \left\lbrace j^*: j^* \in [0 , k] , j^* r - \sum_{i=1}^{j^*} D_{i} + \delta_{j^* , k} I^r_{1,\infty} = W^r_k \right\rbrace.
$$
In words, $i^r_k$ is the (largest) index at which the random walk $W^r_k$
attains its maximum.

From Lemma~\ref{explic1} and the fact that $\{ D_i , i \geq 1\}$ is a sequence of i.i.d. r.v.s,
we obtain the following explicit characterization for the cost incurred by $\pi_r$ over any $L$ consecutive periods.
\begin{lemma} \label{constcost}
{\bf} \hspace*{0.1in}
For any $r \in [0, \E[D])$ and $\tau \in [L+1 , T + 1]$,
\begin{equation*}
\E\left[\sum_{t=\tau}^{\tau + L - 1} C^{\pi_r}_t\right] =
h \sum_{k=1}^{L} \E\left[ i^r_k r - \sum_{i=1}^{i^r_k} D_{i} + \delta_{ i^r_k , k }
I^r_{1,\infty} \right] +  c \big( L \E[D] - L r \big).
\end{equation*}
\end{lemma}

\section{Lower bound on an optimal policy}\label{lowerbound}
We now derive in this section a lower bound on the cost incurred by an optimal policy during any consecutive $L$ time periods.
First, it will be convenient to review a result of Zipkin \cite{Zipkin08}, which establishes an upper bound on the ordering quantities of a family of optimal policies for
Problem~\eqref{probstate}.

\begin{lemma}[\cite{Zipkin08}]\label{zipreview}{\bf} \hspace*{0.1in}
There exists an optimal policy $\overline{\pi}$ for Problem~\eqref{probstate} that with probability (w.p.) $1$ never orders more than $Q$, i.e., $\mathbf{x}^{\overline{\pi}}_t \leq Q \mathbf{1}$ for all $t \in [L + 1 , T]$.
\end{lemma}

\vspace*{0.1in}
\noindent
{\bf Remark.}
Note that Lemma~\ref{zipreview} implicitly asserts the existence of at least one optimal policy for Problem~\eqref{probstate}.
A priori, it could have been possible that no such policy existed, i.e., the optimal value was only attained in the limit by some sequence of policies.
However, the dynamic programming formulation for Problem~\eqref{probstate} provided in \cite{Zipkin08}, combined with the convexity and monotonicity results proven in the same paper, as well as the strict positivity of all underlying cost parameters and continuity of all relevant cost-to-go functions, indeed ensures the existence of at least one optimal policy; we refer the interested reader to \cite{Zipkin08} for details.

\vspace*{0.1in}

We now construct a lower bound by computing the cost that $\overline{\pi}$ incurs over any $L$ consecutive time periods, if the policy were able to choose the state of the
system at the start of those $L$ time periods to be as favorable as possible, subject only to the conditions imposed by Lemma~\ref{zipreview} (which must hold w.p.$1$).

In particular, let $(\mathbf{x}^*, {\mathcal I}^*)$ denote any solution to the optimization problem
\begin{equation}\label{getlow1}
\min_{ \mathbf{x} \in [\mathbf{0}, Q \mathbf{1}] , {\mathcal I} \in \reals^+ }
\E\left[ \sum_{t=\tau}^{\tau + L - 1} C^{\pi}_t \Big| \mathbf{x}^{\pi}_{\tau} =
\mathbf{x}, I^{\pi}_{\tau} = {\mathcal I} \right],
\end{equation}
where the existence of $(\mathbf{x}^*, {\mathcal I}^*)$ follows from the fact that
$\E[ \sum_{t=\tau}^{\tau + L - 1} C^{\pi}_t | \mathbf{x}^{\pi}_{\tau}
= \mathbf{x}, I^{\pi}_{\tau} = {\mathcal I} ]$
is continuous with respect to $(\mathbf{x}, {\mathcal I})$ and goes to infinity as ${\mathcal I}$ goes to infinity, combined with a routine compactness argument.

Note that w.l.o.g.\ we can take $(\mathbf{x}^*, {\mathcal I}^*)$
to be independent of the particular value of $\tau$, and thus a function of ${\mathcal D}, c, h, L$ only.
Further note that conditional on the event $\lbrace \mathbf{x}^{\overline{\pi}}_1 = \mathbf{x}^*, I^{\overline{\pi}}_1 =  {\mathcal I}^* \rbrace$,
the conditional joint distribution of
$\lbrace I^{\overline{\pi}}_{t+1} , N^{\overline{\pi}}_{t} , C^{\overline{\pi}}_{t} ; \, t = 1,\ldots, L \rbrace$
does not depend on the particular policy choices of $\overline{\pi}$, and we denote these conditional r.v.s as $\lbrace I^*_{t+1}, N^*_{t} , C^*_{t} ; \, t = 1,\ldots, L \rbrace$ where $I^*_1 = {\mathcal I}^*$.
Let $$V_k \stackrel{\Delta}{=} \max_{j=0,\ldots,k} \bigg( \sum_{i=k + 1 - j}^k {x}^*_i - \sum_{i= 1}^j D_{i} + \delta_{j,k} {\mathcal I}^* \bigg)$$
and
$$
v^*_k  \stackrel{\Delta}{=}
\max \left\lbrace j^*: j^* \in [0 , k] , \sum_{i=k + 1 - j}^k {x}^*_i - \sum_{i= 1}^j D_{i} + \delta_{j,k} {\mathcal I}^* = V_k \right\rbrace.
$$
Then, combining Lemma~\ref{explic1} with the fact that $\{ D_i, \, i \geq 1\}$ is a sequence of i.i.d. r.v.s, we derive the following lower bound for $\overline{\pi}$.

\begin{lemma}\label{lowerblemma}{\bf} \hspace*{0.1in}
For any $\tau \in [1 , T + 1]$,
\begin{eqnarray*}
\E\left[\sum_{t=\tau}^{\tau + L - 1} C^{\overline{\pi}}_t\right] &\geq&
h \sum_{k=1}^{L} \E[V_k] +  c \big( \E[ I^{*}_{L+1} ] - {\mathcal I}^* + L \E[D] - \sum_{i=1}^{L} x^*_i  \big) .
\end{eqnarray*}
\end{lemma}

Ultimately, we wish to show that the performance of an appropriate constant-order policy nearly matches the lower bound established in Lemma~\ref{lowerblemma}.
To accomplish this, we begin by making several observations.
First, note that if we wanted to select a constant-order policy which came close to matching the aforementioned lower bound, a natural approach would be to select the constant-order policy which ``best mimics'' the pipeline vector $\mathbf{x}^*$, i.e., by considering the policy $\pi_{r^*}$ (assuming $r^* < \E[D]$) and recalling that $r^* = L^{-1} \sum_{i=1}^L x^*_i$.
In this case, we are left with the task of showing that the cost incurred over any consecutive $L$ periods by policy $\pi_{r^*}$, as detailed in Lemma~\ref{constcost}, is ``close'' to the lower bound established in Lemma~\ref{lowerblemma}.

To this end, first note that $\sum_{i=1}^{L} x^*_i = L r^*$.
As such, when comparing the terms appearing in Lemma~\ref{constcost} (applied with constant $r^*$) to those appearing in Lemma~\ref{lowerblemma}, the primary difficulty lies in comparing $\sum_{k=1}^{L} \E[V_k]$ to $\sum_{k=1}^{L} \E[ i^{r^*}_k r^* - \sum_{i=1}^{i^{r^*}_k} D_{i}]$.
We will accomplish this through a particular coupling as follows.
For $k \in [1,L]$, let us construct $V_k, W^{r^*}_k$, and $i^{r^*}_k$ on the same probability space, using the same sequence of demands $\lbrace D_i ; \, i = 1,\ldots,L \rbrace$ (independent of $I^{r^*}_{1,\infty}$).
Since the maximum of several terms is at least any one of the terms (even if selected randomly in an arbitrary manner), for $k \in [1,L]$,
it follows that w.p.$1$
$$V_k \;\; \geq \;\; \sum_{i= k + 1 - i^{r^*}_k }^k x^*_i - \sum_{i = 1}^{i^{r^*}_k} D_{i} + \delta_{i^{r^*}_k , k} {\mathcal I}^*.$$
Upon combining the above with Lemma \ref{lowerblemma} and the non-negativity of $\delta_{i^{r^*}_k , k} {\mathcal I}^*$, we conclude that
\begin{equation}\label{coupleme1}
\sum_{k=1}^{L} \E\left[ i^{r^*}_k r^* - \sum_{i=1}^{i^{r^*}_k} D_{i}\right] - \sum_{k=1}^L \E[V_k] \;\; \leq \;\;
r^* \sum_{k=1}^{L} \E[ i^{r^*}_k] - \sum_{k=1}^L \E[ \sum_{i= k + 1 - i^{r^*}_k }^k x^*_i].
\end{equation}

Note that the intuition behind why the above coupling works is that the indices $j \in [0,k]$ for which $\sum_{i= 1}^j D_{i}$ is exceptionally small are good candidates for both $i^{r^*}_k$ and $v^*_k$.
From a purely technical perspective, the coupling is convenient for two fundamental reasons.
First, it eliminates all terms of the form $\sum_{i= 1}^Z D_{i}$, where $Z$ is a random index that may depend in a complicated way on $\lbrace D_i ; \, i \geq 1 \rbrace$, and is not in general adapted to the filtration generated by $\lbrace D_i ; \, i \geq 1 \rbrace$.
Second, as we shall formalize below, the terms $r^* \sum_{k=1}^{L} \E[ i^{r^*}_k]$ and $\sum_{k=1}^L \E[ \sum_{i= k + 1 - i^{r^*}_k }^k x^*_i]$ can each be well-approximated (in an appropriate sense) by $E[i^{r^*}_{\infty}]  \sum_{i=1}^L x^*_i$.

Combining Lemma~\ref{lowerblemma} with (\ref{coupleme1}), we conclude the following refined lower bound, which will be convenient for comparing the cost incurred by policy $\pi_{r^*}$ to that incurred by $\overline{\pi}$.

\begin{lemma}\label{lowerblemma2}{\bf} \hspace*{0.1in}
For any $\tau \in [1 , T + 1]$, $\E\left[\sum_{t=\tau}^{\tau + L - 1} C^{\overline{\pi}}_t\right]$ is at least
\begin{eqnarray*}
\ &\ &\ h \bigg( \sum_{k=1}^{L} \E\big[ i^{r^*}_k r^* - \sum_{i=1}^{i^{r^*}_k} D_{i}\big] 
+ \sum_{k=1}^L \E[ \sum_{i= k + 1 - i^{r^*}_k }^k x^*_i] - r^* \sum_{k=1}^{L} \E[ i^{r^*}_k] \bigg)
\\&\ &\ \ \ \ \ \ + \ \ \ c \bigg( \E[ I^{*}_{L+1} ] - {\mathcal I}^* + L \E[D] - \sum_{i=1}^{L} x^*_i  \bigg) .
\end{eqnarray*}
\end{lemma}

\section{Difference between constant-order policy and lower bound}\label{boundthedifference}

In this section, we combine Lemmas~\ref{constcost} and \ref{lowerblemma2} to show that the cost incurred by $\pi_{r^*}$ over any $L$ consecutive time periods nearly matches the lower bound established in Lemma~\ref{lowerblemma2}.
We accomplish this through a term-by-term comparison of the expressions appearing in Lemmas~\ref{constcost} and \ref{lowerblemma2}.

\vspace*{0.1in}
\noindent
{\bf Remark.}
We note that conceptually, our approach is closely related to several results in the queueing literature which prove that for certain queueing systems, the arrival process in which all inter-arrival times are the same constant is asymptotically extremal with regards to mean waiting time \cite{hajek1983proof,humblet1982determinism}.  It is an interesting open question to further 
quantify the connection between our approach and the convexity-type arguments typically used to prove extremality in such queueing systems \cite{humblet1982determinism}, which could likely be used to provide an alternate proof of our main results. \vspace*{0.1in}

For $r \in [0,E[D])$, let us define
\begin{equation}
\Theta_r \stackrel{\Delta}{=} \frac{(\E[D]-r)^2}{4(\E^2[D]+ \E[D^2])} .
\label{def:Theta}
\end{equation}
Then the main result of this section is formally stated as follows.

\begin{theorem} \label{afewgoodbounds}
{\bf} \hspace*{0.1in}
If $r^{*} < \E[D]$, then for any $\tau \in [L+1 , T+1]$,
\begin{equation}\label{showsmall1}
\E[\sum_{t=\tau}^{\tau + L - 1} C^{\pi_{r^{*}}}_t] - \E[\sum_{t=\tau}^{\tau + L - 1} C^{\overline{\pi}}_t]
\end{equation}
is at most
$$
h (Q + 2^{\frac{3}{2}} \E[D]) \Theta_{r^*}^{-3} + c {\mathcal I}^*.
$$
\end{theorem}

Before proceeding with the proof of Theorem~\ref{afewgoodbounds}, it will be convenient to derive several bounds for the distribution of $i^r_k$ and $i^r_{\infty}$, where we defer the associated proofs to the technical appendix in Section~\ref{appsec}.

\begin{lemma}\label{boundindex}
{\bf} \hspace*{0.1in}
For any $r \in [0, \E[D])$ and integers $j,k \geq 0$, $i^r_k$ has the same distribution as $\min(k,i^r_{\infty})$, namely $\p(i^r_k = j) = \p\big( \min(k, i^r_{\infty}) = j \big)$.
Furthermore,
$$
\p(i^r_{\infty} \geq  k) \; \leq \; \Theta_r^{-1} \big( 1 - \Theta_r \big)^k \, , \qquad
\sum_{k=0}^{\infty} \sum_{j=k}^{\infty} \p(i^r_{\infty} \geq  j) \; \leq \; \Theta_r^{-3} \, , \qquad
\E[(I^r_{\infty})^2] \; \leq \; 2 \Theta_r^{-3} \E^2[D] .
$$
\end{lemma}

\vspace*{0.1in}
\noindent
{\bf Remark.}
We note that a more precise analysis of the quantities in
Lemma~\ref{boundindex} would be possible using the theory of ladder heights
and epochs~\cite{Asmu03}, especially the precise results for the
relevant moments given in~\cite{Spar53,Spar54,Fell71,Loul78} and the recent
work by Nagaev \cite{Nagaev10R}.
However, since the increments of the random walks that we consider have a very
special structure (i.e., they are absolutely bounded from above), as well as
for the sake of simplicity, we do not pursue such an analysis here.
\vspace*{0.1in}

We now complete the proof of Theorem~\ref{afewgoodbounds}.

\begin{Proof}[Proof of Theorem~\ref{afewgoodbounds}]
Suppose $r^* < \E[D]$.
Combining Lemmas~\ref{constcost} and \ref{lowerblemma2} with the definition of $r^*$ and the non-negativity of all relevant terms,
and then simplifying, allows us to conclude that (\ref{showsmall1}) is at most 
\begin{equation}\label{showsmall2}
h\bigg( r^* \sum_{k=1}^{L} \E[ i^{r^*}_k] - \sum_{k=1}^L \E[ \sum_{i= k + 1 - i^{r^*}_k }^k x^*_i ]
+ \sum_{k=1}^L \E[\delta_{ i^{r^*}_k , k } I^{r^*}_{1,\infty}] \bigg)
+  c {\mathcal I}^*.
\end{equation}

We proceed by bounding each term appearing in (\ref{showsmall2}), beginning with 
\begin{equation}\label{boundfirst1}
 r^* \sum_{k=1}^{L} \E[i^{r^*}_k] - \sum_{k=1}^L \E\bigg[ \sum_{i= k + 1 - i^{r^*}_k }^k x^*_i \bigg].
\end{equation}
First, it will be convenient to generalize our notation $\delta_{j,k}$ as follows.  For an integer $j$ and set $S$, define $\delta_{j,S}$ to be 1 if $j \in S$, and $0$ otherwise.  Observe that
$$
\sum_{i=k + 1 - i^{r^*}_k }^k {x}^{*}_{i}
\;\; = \;\; \sum_{i=1}^k {x}^{*}_{i}  \delta_{i,[k + 1 - i^{r^*}_k, k]}
\;\; = \;\; \sum_{i=1}^k {x}^{*}_{i}  \delta_{i^{r^*}_k , [k+1-i,\infty)},
$$
and thus by interchanging the order of summation and applying Lemma~\ref{boundindex}, we obtain
\begin{eqnarray}
\sum_{k=1}^{L} \E\bigg[\sum_{i=k + 1 - i^{r^*}_k }^k {x}^{*}_{i}\bigg] &=& \sum_{i=1}^L {x}^{*}_{i} \sum_{k=i}^{L} \E[\delta_{i^{r^*}_k , [k+1-i,\infty)}] \nonumber
\\&=&\sum_{i=1}^L {x}^{*}_{i} \sum_{k=i}^{L} \p(i^{r^*}_k  \geq k+1-i)  \nonumber
\\&=& \sum_{i=1}^L {x}^{*}_{i} \sum_{k=1}^{L+1 - i}  \p( i^{r^*}_{k + i - 1} \geq k ) \nonumber
\\&=& \sum_{i=1}^L {x}^{*}_{i} \sum_{k=1}^{L+1 - i} \p\big( i^{r^*}_{\infty} \geq k \big) \label{showsmall3}.
\end{eqnarray}
Applying the definition of $r^*$ together with Lemma~\ref{zipreview} and the fact that $\E[i^{r^*}_k] \leq \E[i^{r^*}_{\infty}]$ by Lemma~\ref{boundindex}, yields that (\ref{boundfirst1}) is at most
\begin{eqnarray}
r^* \sum_{k=1}^{L} \E[i^{r^*}_{\infty}] - \sum_{i=1}^L {x}^{*}_{i} \sum_{k=1}^{L+1 - i} \p\big( i^{r^*}_{\infty} \geq k \big)
&=& \sum_{i=1}^L x^*_i  \sum_{k=1}^{\infty} \p\big( i^{r^*}_{\infty} \geq k \big) - \sum_{i=1}^L {x}^{*}_{i} \sum_{k=1}^{L+1 - i} \p\big( i^{r^*}_{\infty} \geq k \big) \nonumber
\\&=& \sum_{i=1}^L x^*_i  \sum_{k = L + 2 - i}^{\infty} \p\big( i^{r^*}_{\infty} \geq k \big) \nonumber
\\&\leq& Q \sum_{i=1}^L \sum_{k = L + 2 - i}^{\infty} \p\big( i^{r^*}_{\infty} \geq k \big) \nonumber
\\&\leq& Q \sum_{i=0}^{\infty} \sum_{k=i}^{\infty} \p( i^{r^{*}}_{\infty} \geq k ) \nonumber,
\end{eqnarray}
where the final inequality follows from a straightforward reindexing.
Upon combining this with Lemma~\ref{boundindex}, we conclude that (\ref{boundfirst1}) is at most
\begin{equation}\label{boundfirst1done}
Q \Theta_{r^*}^{-3}.
\end{equation}

Next we turn to bound 
\begin{equation}\label{boundsecond1}
\sum_{k=1}^L \E[\delta_{ i^{r^*}_k , k } I^{r^*}_{1,\infty}].
\end{equation}
Applying the Cauchy-Schwartz inequality and Lemma~\ref{boundindex}, yields
\begin{eqnarray}
\sum_{k=1}^L \E[\delta_{ i^{r^*}_k , k } I^{r^*}_{1,\infty}] &\leq& \sum_{k=1}^L \E^{\frac{1}{2}}[\delta_{ i^{r^{*}}_k,k}]\E^{\frac{1}{2}}[(I^{r^{*}}_{\infty})^2 ] \nonumber
\\&=& \sum_{k=1}^L \p^{\frac{1}{2}}\big( i^{r^{*}}_{\infty} \geq k \big) \E^{\frac{1}{2}}[(I^{r^{*}}_{\infty})^2 ] \nonumber
\\&\leq& \sum_{k=1}^{\infty} \big( \Theta_{r^*}^{-1} \big( 1 - \Theta_{r^*} \big)^k \big)^{\frac{1}{2}} \big( 2 \Theta_r^{-3} \E^2[D] \big)^{\frac{1}{2}} \nonumber
\\&=& 2^{\frac{1}{2}} \E[D] \Theta_{r^*}^{-2} \frac{ (1 - \Theta)^{\frac{1}{2}} }{1 - (1 - \Theta)^{\frac{1}{2}}}\ \ \ \leq\ \ \ 2^{\frac{3}{2}} \E[D] \Theta_{r^*}^{-3} \label{showsmall5},
\end{eqnarray}
where the final inequality follows from multiplying and dividing by $1 + (1 - \Theta_{r^*})^{\frac{1}{2}}$ and from noting that $(1 - \Theta_{r^*})^{\frac{1}{2}}\big(1 + (1 - \Theta_{r^*})^{\frac{1}{2}}\big) \leq 2$.

Finally, using both (\ref{boundfirst1done}) to bound (\ref{boundfirst1}) and (\ref{showsmall5}) to bound (\ref{boundsecond1}) in (\ref{showsmall2}), completes the proof.
\end{Proof}

\section{Proof of main result}
\label{provesec}
We now complete the proof of our main result, namely Theorem~\ref{mainresult}, by combining Theorem~\ref{afewgoodbounds} with several additional bounds.
In light of Theorem~\ref{afewgoodbounds}, the primary difficulty which remains is proving that $r^*$ is bounded away from $\E[D]$ as $L \rightarrow \infty$.
Recall that $z = \argmin_{v \geq 0} \left(h \E\big[I^v_{\infty}\big] - c v\right)$.
The final step will be to show that $\pi_{r^*}$ is itself dominated by the policy $\pi_z$, which we will prove to be the ``best-possible'' constant-order policy.
Combining these results with a few straightforward asymptotic arguments will complete the proof.

\subsection{Bounding $r^*$ away from $\E[D]$.}
Let us begin by proving that $r^*$ is bounded away from $\E[D]$ as $L \rightarrow \infty$.
In particular, we will prove the following result, recalling that $m = \Big\lceil \big(26 ( 3 \zeta + c (h \sigma)^{-1} \E[D] + 1 ) \big)^2 \Big\rceil$.

\begin{lemma}\label{boundr}{\bf} \hspace*{0.1in}
For all $L \geq 8 (\sigma^{-1} Q + 1) m^{\frac{3}{2}}$, we have that $\E[D] - r^* \geq \frac{1}{2} \sigma m^{-\frac{1}{2}}$.
\end{lemma}

Although at first glance one might think such a result to be straightforward, the fact that in principle the components of $\mathbf{x}^*$ could vary considerably over the $L$ periods creates difficulties here.
Indeed, we will actually argue indirectly as follows.
Roughly speaking, we argue that for any time $t$, if in the $m$ periods leading up to and including $t$ (i.e., periods $t - m + 1, \ldots, t$) the corresponding pipeline vector  components (i.e., $x^*_{t-m+1},\ldots,x^*_t$) were much larger (on average) than $\E[D]$ (i.e., $\sum_{k = t - m + 1}^t x^*_k - m \E[D]$ is too large), then the expected value of the inventory at the end of period $t$, $I^*_{t+1}$, will be very large.

More precisely, defining $\gamma_t \stackrel{\Delta}{=} \sum_{k = t - m + 1}^t x^*_k$, 
we note that $\E[I^*_{t+1}] \geq \E[\max(0, \gamma_t - \sum_{k = t - m + 1}^t D_k)]$.
We then apply a type of central limit theorem scaling, combined with certain explicit bounds on the rate of convergence in the central limit theorem (i.e., Stein's method), to argue that if $\gamma_t - m \E[D]$ is large, then $\E[I^*_{t+1}]$ is large.
It will follow that $\gamma_t - m \E[D]$ cannot be large for too many values of $t$.
A key insight here is that because of the maximum within the expectation, we get a strong lower bound even when $\gamma_t - m \E[D] = 0$, which in turn allows us to show that $\gamma_t - m \E[D]$ should typically be significantly less than zero.
We then argue that $\sum_{t = m}^L \gamma_t$ is sufficiently close to  $m \sum_{k=1}^L x^*_k$, and combine the above observations to conclude that $\E[D] - r*$ must be bounded away from zero.

\begin{Proof}[Proof of Lemma~\ref{boundr}]
It follows from (\ref{inventory1}) and non-negativity that, for all $t \in [m,L]$, w.p.$1$
$$I^*_t + x^*_t - D_t \;\; \geq \;\; \gamma_t - \sum_{k = t - m + 1}^t D_k.$$
Thus, for all $t \in[m, L]$, $\E[\max(0, I^*_t + x^*_t - D_t)]$ is at least
\begin{eqnarray}
& & \E\bigg[\max\Big(0, \gamma_t - \sum_{k = t - m + 1}^t D_k\Big)\bigg] \nonumber \\
& & \qquad = \quad \sigma m^{\frac{1}{2}} \E\bigg[\max\Big(0, \frac{\sum_{k = t - m + 1}^t (\E[D] - D_k)}{\sigma m^{\frac{1}{2}}} + \frac{\gamma_t - m \E[D]}{\sigma m^{\frac{1}{2}}} \Big)\bigg].
\label{inducto2}
\end{eqnarray}

Let $N$ denote a standard normal r.v.
We now show that (\ref{inducto2}) is well-approximated by 
$$\sigma m^{\frac{1}{2}} \, \E\left[\max\Big(0, \; N + \frac{\gamma_t - m \E[D]}{\sigma m^{\frac{1}{2}}} \Big)\right],$$
using known results on the rate of convergence in the central limit theorem.
Such results are typically derived via Stein's method, and we refer the interested reader to \cite{chen2005stein} for details.
Specifically, the following explicit bound on the rate of convergence in the central limit theorem is generally well known.

\begin{theorem}[\cite{chen2005stein}]\label{Steiner1}{\bf} \hspace*{0.1in}
Suppose that $F: {\mathcal R} \rightarrow {\mathcal R}$ is any Lipschitz-continuous function with Lipschitz constant at most unity, i.e., for all $x,y \in {\mathcal R}$, $|F(x) - F(y)| \leq |x - y|$.
Suppose that $\lbrace X_i; \, i \geq 1 \rbrace$ is any sequence of i.i.d.\ r.v.s s.t.\ $E[X_1] = 0$, $E[X_1^2] = 1$, and $E[|X_1^3|] < \infty$.
Then, for all $n \geq 1$,
$$\bigg|E\Big[F( n^{-\frac{1}{2}} \sum_{i=1}^n X_i)\Big] - E[F(N)] \bigg| \;\; \leq \;\; 3 n^{-\frac{1}{2}} \E[|X_1|^3].$$
\end{theorem}

Letting $F_t(x) \stackrel{\Delta}{=} \max\Big(0, x + \frac{\gamma_t - m \E[D]}{\sigma m^{\frac{1}{2}}} \Big)$, it follows from Theorem~\ref{Steiner1} and (\ref{inducto2}) that, for all $t \in [m, L]$, 
$$\E[\max(0, I^*_t + x^*_t - D_t)] \;\; \geq \;\; \sigma m^{\frac{1}{2}} \left( \E\bigg[\max\Big(0, N + \frac{\gamma_t - m \E[D]}{\sigma m^{\frac{1}{2}}} \Big)\bigg]- 3 m^{-\frac{1}{2}} \zeta \right),$$
and thus
\begin{equation}\label{findit0}
\sum_{t=1}^L \E[C^*_t] \;\; \geq \;\; h \sigma m^{\frac{1}{2}} \sum_{t=m}^L  \E\left[\max\Big(0, N + \frac{\gamma_t - m \E[D]}{\sigma m^{\frac{1}{2}}} \Big)\right] - 3 h \sigma \zeta L.
\end{equation}
Note that $\psi(y) \stackrel{\Delta}{=} \E[\max(0, N + y)]$ is a convex function of $y$.
It then follows from Jensen's inequality that
$$\sum_{t=m}^L  (L - m + 1)^{-1} \E\bigg[\max\Big(0, N + \frac{\gamma_t - m \E[D]}{\sigma m^{\frac{1}{2}}} \Big)\bigg]$$ is at least 
\begin{equation}\label{isbigger1}
\E\left[\max\Big(0, N + (L - m + 1)^{-1} \sum_{t=m}^L \frac{\gamma_t - m \E[D]}{\sigma m^{\frac{1}{2}}} \Big)\right].
\end{equation}
Since $(\mathbf{0},0)$ is a feasible solution to Problem~\eqref{getlow1}, with value $c L \E[D]$, the optimality of $(\mathbf{x}^*, {\mathcal I}^*)$ implies
\begin{equation}\label{beats0}
\sum_{t=1}^L \E[C^*_t] \;\; \leq \;\; L c \E[D].
\end{equation}
Combining (\ref{findit0}), (\ref{isbigger1}), and (\ref{beats0}) yields
\begin{equation}\label{closeto1}
\E\bigg[\max\Big(0, N + (L - m + 1)^{-1} \sum_{t=m}^L \frac{\gamma_t - m \E[D]}{\sigma m^{\frac{1}{2}}} \Big)\bigg] \;\; \leq \;\; \frac{L}{L - m + 1} \big(3 \zeta + c (h \sigma)^{-1} \E[D] \big) m^{-\frac{1}{2}}.
\end{equation}

We now relate $\sum_{t=m}^L \gamma_t$ to $\sum_{t=1}^L x^*_t$ by proving that
\begin{equation}\label{sumofsum1}
\sum_{t = m}^L \gamma_t \;\; = \;\; \sum_{t=m}^L \sum_{k = t - m + 1}^t x^*_k \;\; \geq \;\; m \sum_{k=1}^L x^*_k - 2 m^2 Q.
\end{equation}
Indeed, it follows from a straightforward counting argument that for all $t \in [m, L - m]$, $x^*_t$ appears exactly $m$ times in the double sum $\sum_{t=m}^L \sum_{k = t - m + 1}^t x^*_k$, and thus
$$\sum_{t = m}^L \gamma_t \;\; \geq \;\; m \sum_{t = m}^{L - m} x^*_t.$$
Moreover, since $x^*_t \leq Q$ for all $t$, we conclude that
$$m \sum_{t = m}^{L - m} x^*_t \;\; \geq \;\; m \sum_{k = 1}^L x^*_k - 2 m^2 Q.$$
In combination, these results yield (\ref{sumofsum1}).

Next, upon combining (\ref{closeto1}) and (\ref{sumofsum1}) with the monotonicity of $\psi$, and simplifying all relevant expressions, we obtain
\begin{equation}\label{closeto2}
\E\left[\max\bigg(0, N + \frac{m^{\frac{1}{2}}}{\sigma} \Big( \frac{L}{L - m  +1} r^* - \E[D] - \frac{2 m Q}{L - m +1} \Big) \bigg)\right] \; \leq \; \frac{L}{L - m + 1} \big(3 \zeta + c (h \sigma)^{-1} \E[D] \big) m^{-\frac{1}{2}}.
\end{equation}
Noting that $L \geq 2 m$ implies $\frac{L}{L - m + 1} \leq 2$ and $L \geq 4 \sigma^{-1} m^{\frac{3}{2}} Q$ implies $\frac{2 m Q}{L} \leq \frac{1}{2} \sigma m^{-\frac{1}{2}}$, we devise from the monotonicity of $\psi$ that, for all $L \geq \max(2 m, 8 \sigma^{-1} m^{\frac{3}{2}} Q)$, 
\begin{equation}\label{closeto3}
\E\left[\max\bigg(0, N + \frac{m^{\frac{1}{2}}}{\sigma} \Big( r^* - \E[D] - \frac{1}{2} \sigma m^{-\frac{1}{2}} \Big) \bigg)\right] \;\; \leq \;\; 2 \big(3 \zeta + c (h \sigma)^{-1} \E[D] \big) m^{-\frac{1}{2}}.
\end{equation}
It is easily verified that $(\E[\max(0, N - 1)])^{-1} \leq 13$, and thus from definitions and basic algebra
$$m \;\; \geq \;\; \bigg(2 \big(3 \zeta + c (h \sigma)^{-1} \E[D] \big) (\E[\max(0, N - 1)])^{-1}\bigg)^2.$$
We conclude that the right-hand-side of (\ref{closeto3}) is at most $\E[\max(0, N - 1)]$, in which case the monotonicity of $\psi$ implies
$$\frac{m^{\frac{1}{2}}}{\sigma} \left( r^* - \E[D] - \frac{1}{2} \sigma m^{-\frac{1}{2}} \right) \;\; \leq \;\; - 1,$$
namely $\E[D] - r^* \geq \frac{1}{2} \sigma m^{-\frac{1}{2}}$.
Combining the above with some straightforward algebra completes the proof.
\end{Proof}

We end this subsection by combining Theorem~\ref{afewgoodbounds} and Lemma~\ref{boundr} to bound the difference between the constant-order policy and an optimal policy.
\begin{corollary}\label{endit1}
{\bf} \hspace*{0.1in}
For all $L \geq 8 (\sigma^{-1} Q + 1) m^{\frac{3}{2}}$ and any $\tau \in [L+1 , T+1]$, (\ref{showsmall1}) is at most
$$2^{12} h (Q + 2^{\frac{3}{2}} \E[D])(\E^2[D] + \E[D^2])^3 \sigma^{-6} m^3 + c  \big(\lceil (2 c h^{-1} L )^{\frac{1}{2}} \rceil + 2 \big).$$
\end{corollary}
\begin{Proof}
It follows from Lemma~\ref{boundr} that
$$h (Q + 2^{\frac{3}{2}} \E[D]) \Theta_{r^*}^{-3} \;\; \leq \;\; 2^{12} h (Q + 2^{\frac{3}{2}} \E[D])(\E^2[D] + \E[D^2])^3 \sigma^{-6} m^3.$$
Thus, by Theorem~\ref{afewgoodbounds}, to complete the proof it suffices to demonstrate that ${\mathcal I}^* \leq \lceil (2 c h^{-1} L )^{\frac{1}{2}} \rceil + 2$.
Indeed, suppose for contradiction that
${\mathcal I}^* > \big( \lceil (2 c h^{-1} L)^{\frac{1}{2}} \rceil + 2 \big) \E[D].$
For all $k\in [1, \lceil (2 c h^{-1} L )^{\frac{1}{2}} \rceil + 2]$, we have
$$\E[I^*_{1 + k}] \ge {\mathcal I}^* - k \E[D]> \E[D] \big( \lceil (2 c h^{-1} L)^{\frac{1}{2}} \rceil + 2 -  k \big).$$
The resulting holding costs ensure that
$\sum_{t= 1}^L \E[ C^*_t]$ is strictly greater than
\begin{eqnarray*}
h  \E[D] \sum_{k=1}^{\lceil (2 c h^{-1} L)^{\frac{1}{2}} \rceil} k
&\geq& c \E[D] L.
\end{eqnarray*}
Combining this with (\ref{beats0}) completes the proof.
\end{Proof}

\subsection{Proof of Theorem~\ref{mainresult}.}
With Corollary~\ref{endit1} in hand, we now proceed with the proof of our main result, i.e., Theorem~\ref{mainresult}.

\begin{Proof}[Proof of Theorem~\ref{mainresult}]
Suppose $T \geq L$, $\epsilon \in (0,1)$, and $L \geq 8 (\sigma^{-1} Q + 1) m^{\frac{3}{2}}$.
It then follows from Lemma~\ref{boundr} that $r^* < \E[D]$.
Note that 
\begin{equation}\label{tobound1}
\frac{\sum_{t = L + 1}^{T + L} \E[ C^{\pi_{r^*}}_t]}{\sum_{t = L + 1}^{T + L} \E[ C^{\overline{\pi}}_t]}
\end{equation}
equals
$$
\frac{ \sum_{k=1}^{\lfloor \frac{T}{L} \rfloor} \sum_{t = k L + 1}^{(k+1) L} \E[ C^{\pi_{r^*}}_t] + \sum_{t = (\lfloor \frac{T}{L} \rfloor  + 1) L + 1}^{T + L} \E[ C^{\pi_{r^*}}_t] }{ \sum_{k=1}^{\lfloor \frac{T}{L} \rfloor} \sum_{t = k L + 1}^{(k+1) L} \E[ C^{\overline{\pi}}_t] + \sum_{t = (\lfloor \frac{T}{L} \rfloor  + 1) L + 1}^{T + L} \E[ C^{\overline{\pi}}_t] }.
$$
As the policy $\pi_{r^*}$ is stationary and yields a stationary sequence of inventory and ordering levels, it follows that $\E[ C^{\pi_{r^*}}_t] = \E[ C^{\pi_{r^*}}_{L+1}]$ for all $t \geq L + 1$.
Further note that $\E[C^{\overline{\pi}}_t] \geq g$ for all $t \geq L + 1$, and thus $\sum_{t = k L + 1}^{(k+1) L} \E[ C^{\overline{\pi}}_t] \geq L g$ for all $k \in [1, \lfloor \frac{T}{L} \rfloor]$.
Combining the above with Corollary~\ref{endit1}, and the non-negativity of all relevant terms, 
we conclude that (\ref{tobound1}) is at most
\begin{equation}\label{soclose1}
\frac{ \lceil \frac{T}{L} \rceil }{\lfloor \frac{T}{L} \rfloor} \left( 1 + \Big( 2^{12} h (Q + 2^{\frac{3}{2}} \E[D])(\E^2[D] + \E[D^2])^3 \sigma^{-6} m^3 + c  \big(\lceil (2 c h^{-1} L )^{\frac{1}{2}} \rceil + 2 \big) \Big) ( g L)^{-1} \right).
\end{equation}

Next, we note that
$L \; \geq \; 2^{14} h (Q + 2^{\frac{3}{2}} \E[D])(\E^2[D] + \E[D^2])^3 \sigma^{-6} m^3 g^{-1} \epsilon^{-1}$~~~~implies
$$2^{12} h (Q + 2^{\frac{3}{2}} \E[D])(\E^2[D] + \E[D^2])^3 \sigma^{-6} m^3  ( g L)^{-1} \; \leq \; \frac{\epsilon}{4},$$
$$L \; \geq \; \left( 12 c g^{-1} \big( (2 c h^{-1})^{\frac{1}{2}} + 3 \big) \right)^2 \epsilon^{-2} \qquad \mbox{implies} \qquad
c  \big(\lceil (2 c h^{-1} L )^{\frac{1}{2}} \rceil + 2 \big) ( g L)^{-1} \; \leq \; \frac{\epsilon}{12},$$
and 
$$
T \; \geq \; \left(1 + \frac{3}{\epsilon}\right) L \qquad \mbox{implies} \qquad \frac{ \lceil \frac{T}{L} \rceil }{\lfloor \frac{T}{L} \rfloor} \; \leq \; 1 + \frac{\epsilon}{3}.
$$
Combining the above with the fact that $(1 + \frac{\epsilon}{3})^2 \leq 1+ \epsilon$ for all $\epsilon \in (0,1)$ completes the proof that the stated performance guarantees are attained by the policy $\pi_{r^*}$, for any $L,T$ satisfying the conditions of Theorem~\ref{mainresult}.

The final step is to prove that the same guarantees extend to $\pi_z$.
Indeed,  it follows from stationarity that, for any $r \in [0, \E[D])$ and $t \geq L + 1$,
\begin{equation}\label{anyconstant0}
\E[C^{\pi_r}_t] \;\; = \;\; h \E[(I^r_{\infty} + r - D)^+] + c \E[(I^r_{\infty} + r - D)^{-}].
\end{equation}
However, since $I^r_{\infty}$ is distributed as $(I^r_{\infty} + r - D)^+$, we conclude
$$\E[(I^r_{\infty} + r - D)^+] \;\; = \;\; \E[I^r_{\infty}].$$
Furthermore, since
$$\E[(I^r_{\infty} + r - D)^{+}] - \E[(I^r_{\infty} + r - D)^{-}] \;\; = \;\; \E[I^r_{\infty}] + r - \E[D],$$
it follows that
$$\E[(I^r_{\infty} + r - D)^{-}] \;\; = \;\; \E[D] - r.$$
In combination with (\ref{anyconstant0}), we have that for any $r \in [0,\E[D])$ and $t \geq L + 1$
\begin{equation}\label{anyconstant1}
\E[C^{\pi_r}_t] \;\; = \;\; h \E[I^r_{\infty}] + c \big( \E[D] - r \big).
\end{equation}
The desired result then follows from the fact that $z$ is a minimizer of (\ref{anyconstant1}), completing the proof.
\end{Proof}

\section{Conclusion}
\label{conclusion}
In this paper, we considered the single-item, periodic-review, lost-sales
model with positive lead times and i.i.d.\ demand, for which the optimal policy is poorly understood and computationally intractable.
We proved that, as the lead time grows (with the demand distribution,
lost-sales penalty, and holding cost remaining fixed), a simple, open-loop
constant-order policy is in fact asymptotically optimal.
We also established explicit bounds on how large the lead time should be to
ensure that the best constant-order policy incurs an expected cost of at most
$1+\epsilon$ times that incurred by the optimal policy.
To the best of our knowledge, this is the first algorithm proven to be within
$1 + \epsilon$ of optimal for lost-sales models when the lead time is large,
while maintaining a runtime that does not grow with the lead time.
Our main proof technique involved a novel coupling for suprema of random walks,
and may be useful in other settings.

This work leaves many interesting directions for future research.
We suspect that our explicit bounds are not tight, and a more precise analysis
of the constant-order policy would further help to explain the good performance of the algorithm for lead times as
small as four, as reported by Zipkin \cite{Zipkin08b}.
Since lost sales models commonly arise in practice, an interesting challenge
is to combine the core ideas of our analysis with known results from dynamic
programming to derive and analyze practical ``hybrid'' algorithms, which use
more elaborate forms of dynamic programming when the lead time is small and
gradually transition to less computationally intensive algorithms (with the
constant-order policy at the extreme) as the lead time grows.
It would also be interesting to prove that a similar phenomenon occurs for other policies, as well as in other
inventory models.
In particular, it is an interesting open question whether other simple (but perhaps slightly more sophisticated) policies, such as the order-up-to policy considered by Huh et al. \cite{huh2009asymptotic} and the cost-balancing policy considered by Levi et al. \cite{Levi08c}, exhibit a similar asymptotic optimality as the lead time grows.
On a related note, the fact that $y(\epsilon)$ contains terms of the form $c h^{-1}$ suggests that the constant-order policy may require larger lead times to approach near-optimality when the ratio of lost-sales-penalty to holding-cost is large.
Since this is exactly the regime in which the order-up-to policy of \cite{huh2009asymptotic} provably works well, it is an interesting open question to try and combine these (and perhaps other) algorithms to yield tighter performance guarantees over a larger range of parameters.

Philosophically, our main results and insights fall under the broad heading of ``long-range
independence / decay of correlations'' phenomena, in which so much uncertainty is introduced into a model that even very sophisticated algorithms cannot perform significantly better than very simple algorithms.
Such ideas have led to significant progress on fundamental models in other fields \cite{mossel2012stochastic, massoulie2013community, gamarnik2013correlation, sly2010computational, gamarnik2014limits, deshpande2013finding, feldman2013complexity}, and may prove useful in other operations management problems.

\section{Appendix}\label{appsec}
\subsection{Proof of Lemma~\ref{boundindex}.}
In this appendix, we provide the proof of Lemma~\ref{boundindex}.

\begin{Proof}[Proof of Lemma~\ref{boundindex}.]  We first prove that $i^r_k$ has the same distribution as $\min(k,i^r_{\infty})$.
Let $\lbrace {D'}_i ; \, i \geq 1 \rbrace$ be an additional sequence of i.i.d.\ realizations from ${\mathcal D}$, mutually independent
from $\lbrace D_i ; \, i \geq 1 \rbrace$.
Then, for any $k \geq 1$, we can construct
$I^r_{1,\infty}, W^r_k, \lbrace D_i  ; \, i \geq 1\rbrace, \lbrace D'_i  ; \, i \geq 1\rbrace$ on the same
probability space s.t.\ $I^r_{1,\infty} =  \max_{j \geq 0} ( j r - \sum_{i=1}^j {D'}_i ).$
It is easy to see that we only need to show
$$\p(i^r_k=j) \;\; = \;\; \p(i^r_{\infty} = j) \; , \qquad j = 0,\ldots,k - 1.$$
By definition, we have
$\p[i^r_k=j] \, = \, \p[ {\mathcal I}_1 \cap {\mathcal I}_2\cap {\mathcal I}_3]$
where
\begin{align*}
{\mathcal I}_1 &= \left\{  \ell r - \sum_{i=1}^\ell D_{i} \le  j r - \sum_{i=1}^j D_{i} , \;\; \forall \ell \le j \right\} , \\
{\mathcal I}_2 & = \left\{\ell r - \sum_{i=1}^\ell D_{i} <  j r - \sum_{i=1}^j D_{i}, \;\; \ell =j+1, \ldots, k-1  \right\} , \\
{\mathcal I}_3 & = \left\{ j r - \sum_{i=1}^j D_{i} >  k r - \sum_{i=1}^k D_{i} + I^r_{1,\infty}\right\} .
\end{align*}
Note that
\begin{align*}
{\mathcal I}_3 & = \left\{ j r - \sum_{i=1}^j D_{i} >  k r - \sum_{i=1}^k D_{i} +\max_{\ell \geq 0} \left( \ell r - \sum_{i=1}^\ell {D'}_i  \right)\right\}\\
& =  \left\{j r - \sum_{i=1}^j D_{i} >  \max_{\ell \geq 0}  \left[ (k+\ell) r - \left(\sum_{i=1}^k D_{i}  + \sum_{i=1}^\ell {D'}_i \right) \right]\right\}.
\end{align*}
It therefore follows, since $\lbrace D_i; i \geq 1 \rbrace$ and $\lbrace D'_i; i \geq 1 \rbrace$ are mutually independent i.i.d.\ sequences with common distribution ${\mathcal D}$, that 
$$\p[i^r_k=j] \;\; = \;\; \p[ {\mathcal I}_1 \cap {\mathcal I}_2\cap {\mathcal I'}_3],$$
where
\begin{equation*}
{\mathcal I'}_3 \;\; =  \;\; \left\{j r - \sum_{i=1}^j D_{i} >  \max_{\ell \geq k}  \left[ \ell r - \sum_{i=1}^\ell D_i  \right]\right\} .
\end{equation*}
Noting that this is the definition of $\p(i^r_{\infty}=j)$ completes the first part of the proof.

Before proving the remainder of the lemma, it will be useful to establish that
\begin{equation}\label{provefirst}
\p(i^r_{\infty} = k) \;\; \leq \;\; \big( 1 - \Theta_r \big)^k.
\end{equation}
By definition,
$$\p(i^r_{\infty} = k) \;\; \leq \;\; \p\bigg( \sum_{i=1}^k (r - D_i) \geq 0 \bigg).$$
Applying a Chernoff bound, we find that for any $\theta > 0$
$$\p(i^r_{\infty} = k) \;\; \leq \;\; \E^k\big[\exp\big(\theta(r-D)\big)\big]$$
where
\begin{eqnarray*}
\E\big[\exp\big(\theta(r - D) \big) \big]
& \; = \; & \exp(\theta r)\E\big[\exp\big(- \theta D\big) \big] \\
& \; \leq \; & \exp(\theta r)\E\big[\big(1 + \theta D\big)^{-1} \big] \qquad\qquad (\mbox{since $\exp(v) \geq 1 + v$}) \\
& \; \leq \; & \E\left[\frac{1 + \theta r + \theta^2 r^2 }{1 + \theta D} \right], \qquad\qquad\qquad \mbox{for all $\theta \in (0,r^{-1}]$},
\end{eqnarray*}
the final inequality following from a simple Taylor-series expansion.
However, w.p.$1$, we have
\begin{eqnarray*}
\frac{1 + \theta r + \theta^2 r^2 }{1 + \theta D} &\; =\; & 1 + \theta(r-D) + \frac{\theta^2}{1 + \theta D} \big(r^2 - D(r - D) \big) \\
&\; \leq\; & 1 + \theta(r-D) + \theta^2 (r^2 + D^2),
\end{eqnarray*}
and thus
$$\E\big[\exp\big(\theta(r - D) \big) \big] \;\; \leq  \;\; 1 - \theta\big( \E[D] - r \big) + \theta^2 \big( r^2 + \E[D^2] \big).$$
Observing that
\begin{eqnarray*}
\frac{\E[D]-r}{2(r^2 + \E[D^2])} &\; \leq\; & \frac{\E[D]}{2\E[D^2]} \\
&\; \leq\; & \frac{\E[D]}{2 \E^2[D]} \;\; = \;\; \frac{1}{2 \E[D]} \;\; < \;\; r^{-1},
\end{eqnarray*}
we may take
$$\theta \;\; = \;\; \theta^* \;\; \stackrel{\Delta}{=} \;\; \frac{\E[D]-r}{2(r^2 + \E[D^2])}$$
to conclude
\begin{eqnarray*}
\E\big[\exp\big(\theta^*(r - D) \big) \big] &\; \leq\; & 1 - \frac{ (\E[D] - r)^2 }{4(r^2 + \E[D^2])} \\
&\; \leq\; & 1 - \frac{ (\E[D] - r)^2 }{4(\E^2[D] + \E[D^2])} \;\; = \;\; 1 - \Theta_r,
\end{eqnarray*}
where the final inequality follows from the fact that $r^2 \leq \E^2[D]$.
Combining the above completes the proof of (\ref{provefirst}).

With (\ref{provefirst}) in hand, the lemma follows directly from the basic manipulation of a few geometric series
and the fact that $I^r_{\infty} \leq r i^r_{\infty} \leq \E[D] i^r_{\infty}$ w.p.$1$, the details of which we omit.
\end{Proof}

\section*{Acknowledgments.}
The authors thank Maury Bramson, Jim Dai, Bruce Hajek, Retsef Levi, Marty Reiman, Alan Scheller-Wolf, and Paul Zipkin for several stimulating discussions.
We also thank the anonymous referees and the editors for their constructive feedback on an earlier version of this work.
The first author also thanks the IBM T.J.\ Watson Research Center, for providing a great research environment.

\bibliographystyle{amsplain}
\bibliography{Large_Lead_Times_8_21}





\end{document}